\documentclass[11pt,a4paper,reqno]{article}
\usepackage[T1]{fontenc}
\usepackage[utf8]{inputenc}
\usepackage[english]{babel}

\usepackage{microtype,amsmath,amssymb,amsthm,mathtools}
\usepackage{lmodern}
\usepackage[margin=27mm,headheight=14pt]{geometry}
\usepackage{enumitem,xcolor,fancyhdr}
\usepackage[colorlinks=true,linkcolor=blue!40!black,citecolor=blue!40!black,urlcolor=blue!40!black]{hyperref}
\usepackage{bookmark}

\setlist{itemsep=4pt,topsep=5pt}

\newtheorem{theorem}{Theorem}[section]
\newtheorem{lemma}[theorem]{Lemma}

\theoremstyle{definition}

\theoremstyle{remark}
\newtheorem{remark}[theorem]{Remark}
\numberwithin{equation}{section}

\newcommand{\B}{\mathcal B(X)}

\newcommand{\NN}{\mathcal N(X)}

\newcommand{\Fone}{\mathcal F_1}
\newcommand{\None}{\mathcal N_1}
\newcommand{\Ione}{\mathcal I_1}

\newcommand{\Ker}{\operatorname{Ker}}
\newcommand{\rank}{\operatorname{rank}}
\newcommand{\tr}{\operatorname{tr}}
\newcommand{\C}{\mathbb C}

\title{Surjective Maps Preserving Products with Square Zero}
\author{Khalid Amnchar, Zakariae Aznay, Adel Babbah and Fadil Chabbabi\\[1mm]
\small Laboratory (A.F.N.L.A), Department of Mathematics, Faculty of Sciences,\\
\small Abdelmalek Essa\^adi University, Tetouan, Morocco\\[1mm]
\small \texttt{f.chabbabi@uae.ac.ma}}
\date{}

\begin{document}
\maketitle

\begin{abstract}
Let $X$ be a complex Banach space with $\dim X\geq 3$, and let
$\mathcal{B}(X)$ denote the algebra of bounded linear operators on
$X$. We characterize the surjective maps
$\Phi:\mathcal{B}(X)\to\mathcal{B}(X)$ satisfying
\[
(AB)^2=0
\quad\Longleftrightarrow\quad
(\Phi(A)\Phi(B))^2=0,
\qquad A,B\in\mathcal{B}(X).
\]
More precisely, in the infinite-dimensional case, there exist a
bounded bijective linear or conjugate-linear operator $S:X\to X$
and a function $c:\mathcal{B}(X)\setminus\{0\}\to
\mathbb{C}\setminus\{0\}$ such that
\[
\Phi(A)=c(A)SAS^{-1},\qquad A\neq 0.
\]
In the finite-dimensional case, identifying $\mathcal{B}(X)$ with
$M_d(\mathbb{C})$ and  $d\geq 3$, there exist an invertible matrix
$S\in M_d(\mathbb{C})$, a field automorphism
$\sigma:\mathbb{C}\to\mathbb{C}$, and a function
$c:M_d(\mathbb{C})\to\mathbb{C}\setminus\{0\}$
such that
\[
\Phi(A)=c(A)SA^\sigma S^{-1},\qquad A\neq 0,
\]
where, $A^\sigma=(\sigma(a_{ij}))$ whenever $A=(a_{ij})\in M_d(\C)$.
\end{abstract}


\section{Introduction and main results}
Throughout the paper, $X$ denotes a complex Banach space of dimension greater than $2$, 
$\mathcal{B}(X)$ the algebra of all bounded linear operators on $X$, 
and $X^*$ its continuous dual. When $X$ is finite-dimensional, with 
$\dim X=d$, we identify $X$ with $\mathbb{C}^d$ and 
$\mathcal{B}(X)$ with $M_d(\mathbb{C})$, the algebra of $d\times d$ 
complex matrices.

An operator $A$ is said to be quasi-nilpotent when its spectral radius $r(A)$ is zero, and $A$ is called nilpotent if there exists a non-negative integer $n \in \mathbb{N}$ such that $A^n = 0$. The smallest nonzero integer $n$ satisfying $A^n = 0$ is called the order of nilpotency of $A$. We denote the sets of all quasi-nilpotent and nilpotent operators by $\mathcal{Q}$ and $\NN$, respectively. An operator is said to have square zero if it is nilpotent of order $2$.

The study of maps preserving mathematical relations is a vast and active domain in mathematics (see, e.g.,~\cite{Botta, Li, Li2005, Molnar2002a, molnar2007selected, Semrl2004, Wang2009}). In particular, the preservation of properties and quantities associated with operator or matrix products has attracted considerable attention from numerous researchers over recent years (see, for instance,~\cite{bresar2004mappings, bresar2012multiplication, Bresar, li2012nonsurjective, Semrl,  semrl2008zero}).

In \cite[Theorem 2.1]{Semrl}, Li, \v{S}emrl, and Sze characterized the structure of surjective maps $\phi:\mathcal{B}(X)\to\mathcal{B}(X)$ that preserve the nilpotence of products in the following sense
\begin{equation*}\label{eq:SS}
AB\in \NN\quad\Longleftrightarrow\quad \phi(A)\phi(B)\in \NN
\qquad (A,B\in\mathcal{B}(X)).
\end{equation*}
In  \cite{AC}, the authors give a complete description of of surjective maps $\phi:\mathcal{B}(X)\to\mathcal{B}(X)$, that satisfies 
\begin{equation*}\label{eq:ZF}
AB\in \mathcal{Q}\quad\iff\quad \phi(A)\phi(B)\in  \mathcal Q
\qquad (A,B\in\mathcal{B}(X)).
\end{equation*}
Following this line of research, in the present paper,  we investigate surjective maps that
preserve the square-zero property of products in both directions.
More precisely, we study the structure of surjective maps
$\Phi:\mathcal{B}(X)\to\mathcal{B}(X)$ satisfying
\begin{equation}\label{eq:P}
(AB)^2=0
\quad\Longleftrightarrow\quad
(\Phi(A)\Phi(B))^2=0,
\qquad A,B\in\mathcal{B}(X).\qquad 
\tag{1}
\end{equation}

The only global assumption imposed on $\Phi$ is surjectivity. In
particular, we do not assume additivity, continuity, or injectivity.
The condition in \eqref{eq:P} differs from the preservation of the
nilpotency of $AB$, without specifying its index, as studied in
\cite{Semrl}.

The following two theorems are the main results in this paper. 
\begin{theorem}\label{thm:infinite}
Let $X$ be an infinite-dimensional complex Banach space. A surjective map $\Phi:\B\to\B$ satisfies \eqref{eq:P} if and only if $\Phi(0)=0$ and there exist a bounded invertible linear or conjugate-linear operator $S:X\to X$ and a function
$c:\B\setminus\{0\}\to\C\setminus\{0\}$ such that
\begin{equation}\label{eq:infiniteform}
 \Phi(A)=c(A)SAS^{-1}, \qquad(A\ne0).
\end{equation}
Moreover, for each fixed $A\ne0$,   the map $\lambda\in \C^* \mapsto c(\lambda A)\zeta(\lambda)\in \C^*$ is surjective, where $\zeta(\lambda)=\lambda$ in the case when $S$ is linear and $\zeta(\lambda)=\overline{\lambda}$ when $S$ is conjugate linear, with $\C^*=\C\setminus\{0\}$. 
\end{theorem}
In the finite-dimensional case, we have the following theorem.
\begin{theorem}\label{thm:finite}
Let $d \geq 3$. A surjective map $\Phi : M_d(\mathbb{C}) \to M_d(\mathbb{C})$ satisfies \eqref{eq:P} if and only if there exist an invertible matrix $S \in GL_d(\mathbb{C})$ and a field automorphism $\sigma$ of $\mathbb{C}$ such that, for every nonzero $A \in M_d(\mathbb{C})$,
\[
\Phi(A) = c(A) S A^{\sigma} S^{-1},
\]
where $A^\sigma = (\sigma(a_{ij}))$ for $A = (a_{ij})$, and $c(A) \in \mathbb{C}^*$ is a scalar function satisfying
\begin{equation}\label{eq:scalardim}
 \{c(\lambda A)\sigma(\lambda):\lambda\in\C^*\}=\C^*, 
\end{equation}
\end{theorem}

The paper is organized as follows. 
Section~\ref{section:r1} is devoted to recalling some notation and basic properties of rank-one operators. We also recall a useful relation between rank-one operators that have been introduced in the literature. Furthermore, we provide characterizations of rank-one operators and rank-one nilpotent operators in terms of square-zero products. In particular, Lemma~\ref{lem:1} establishes a necessary condition for two operators to be proportional, expressed in terms of square-zero products involving rank-one operators.  We also introduce the class of operators that can be written as the sum of a scalar multiple of the identity and a rank-one nilpotent operator. We then characterize this class in terms of the following set, which will play an important role throughout the paper:
\[
\Gamma(C):=\{R\in\mathcal{N}_1(X):(CR)^2=0\}.
\]

Section~\ref{section:lemmas} is devoted to establishing several auxiliary results concerning surjective maps \(\Phi:\B\to\B\) satisfying condition~\eqref{eq:P} in both directions. We first show that such maps preserve proportionality between operators, as well as the set of nilpotent operators of order at most four, in both directions; see Lemma~\ref{prop:projective}.

Furthermore, we prove that \(\Phi\) preserves, in both directions, the sets of rank-one operators and rank-one nilpotent operators; see Lemmas~\ref{cor:nilpreserved} and~\ref{thm:rankone}. We then show that, after multiplication by a suitable nonzero scalar-valued function, the resulting map preserves the set of rank-one idempotents and the orthogonality relation between rank-one idempotents.

By applying well-known results from the literature on transformations preserving rank-one idempotents and orthogonality, we obtain a precise description of the action of \(\Phi\) on the set of rank-one idempotents. Finally, in the last section~\ref{section:proof}, we extend this description to the whole algebra \(\B\).

\section{Preliminary Results on Rank-One Operators}\label{section:r1}

For $x\in X$ and $f\in X^*$,  the rank-one operator $x\otimes f$ is defined by
\[
 (x\otimes f)(y)=f(y)x\qquad(y\in X).
\]
Every operator of rank at most one can be represented in this form.  We denote by  $\Fone(X)$, $\None(X)$ and $\Ione(X)$  the sets of rank-one operators, rank-one nilpotent operators, and rank-one idempotents in $\B$, respectively. Thus,
\[
 \None(X)=\{N\in\Fone(X):N^2=0\},\quad\text{ and} \quad \Ione(X)=\{P\in\Fone(X):P^2=P\}.
\]
Recall that, for every operator $T=x\otimes f$ of rank at most one, we have 
\begin{equation}\label{eq:trace}
    T^2=\tr(T)T=f(x)T.  
\end{equation}
Let us recall here a relation on the set of rank-one operators introduced in \cite{Semrl}. 
For two rank-one operators $R_1,R_2\in \Fone(X)$, we write $R_1\sim R_2$ if they have
the same range or the same kernel. In particular, if  $R_1=x\otimes f$ and $R_2=y\otimes g$, then
\(
x\otimes f\sim y\otimes g
\)
if and only if $x$ and $y$ are linearly dependent or $f$ and $g$
are linearly dependent. Moreover, if $R_1\sim R_2$, then $\lambda R_1\sim\beta R_2$
for all nonzero scalars $\lambda,\beta\in\mathbb C$.

For an operator $A\in\B$, we denote by 
\[
 Z(A)=\{B\in\B:(AB)^2=0\}.
\]
In particular, 
for every $A\in\B$ and for every $T=x\otimes f \in \Fone(X)$
\begin{equation}\label{eq:rankoneidentity}
 (AT)^2=(A(x\otimes f))^2=f(Ax)\,Ax\otimes f.
\end{equation}
Therefore,
\[T \in Z(A)\quad\Longleftrightarrow\quad
 (A(x\otimes f))^2=0\quad\Longleftrightarrow\quad f(Ax)=0.
\]
 Consequently, $Z(0)=\B$ and  $Z(A)=\B$ if and only if $A=0$.

The following  lemmas are taken from \cite{Semrl}, they  provides a
necessary and sufficient condition for two rank-one operators
to be in relation.

\begin{lemma}\cite{Semrl}\label{lem:adjacency}
Let $R_1, R_2\in\Fone(X)$ be two linearly independent rank-one operators. Then $R_1\sim R_2$ if and only if there exists $R_3 \in \Fone(X)$ such that $R_3$ is linearly independent of both $R_1$ and $R_2$, and
\[
 Z(R_1)\cap Z(R_2)\subseteq Z(R_3).
\]
\end{lemma}

\begin{lemma}\cite{Semrl}\label{lem:orthogonality}
Let $P,Q\in\Ione(X)$ be distinct. Then $PQ=QP=0$ if and only if there are rank-one nilpotents $M,N$ such that
\[
P\sim M,\quad P\sim N,\quad Q\sim M,\quad Q\sim N,
\qquad M\not\sim N.
\]
\end{lemma}

The following lemma characterizes the set of rank-one idempotent operators in terms of square-zero products.
\begin{lemma}\label{lem:nilpotents}
If $\dim X\ge3$, a nonzero operator $A\in\B$ belongs to $\None(X)$ if and only if
\begin{equation}\label{eq:Nil_one}
 A^4=0\quad\text{and}\quad
 \bigl((AB)^2=0\Longleftrightarrow(BA)^2=0\bigr)
 \quad\text{for every }B\in\B.
\end{equation}
If $\dim X=2$, then $A\in\None(X)$ if and only if $A\ne0$ and $A^4=0$.
\end{lemma}

\begin{proof}
First, suppose that  $A=x\otimes f$ is of  rank one. Hence, for every $B\in \B$, the operators $AB$ and $BA$ have rank at most one, and their traces are both equal to $f(Bx)$. Therefore, by (\ref{eq:trace}), it follows that $(AB)^2=0$ if and only if $(BA)^2=0$. Moreover, if $A\in \None(X) $, then $A^2=0$ and thus $A^4=0$ is also.

Conversely, suppose $A^4=0$ and $\rank A\ge2$. Choose $0\ne x_0\in\Ker A$, and choose $x_1,x_2\in X$ so that $y_1=Ax_1$ and $y_2=Ax_2$ are linearly independent. By Hahn--Banach, there are $f_1,f_2\in X^*$ with $f_i(y_j)=\delta_{ij}$. Put
\[
 B=x_0\otimes f_1+x_1\otimes f_2.
\]
Then, $AB=y_1\otimes f_2$, so $(AB)^2=0$.  On the other hand  $BAx_2=x_1$ and $BAx_1=x_0$, and hence $(BA)^2x_2=x_0\ne0$. This contradicting \eqref{eq:Nil_one}. Thus $\rank A=1$.
In dimension two, every nonzero nilpotent operator has rank one, so the last assertion follows as well.
\end{proof}

The following lemma characterizes when two nonzero operators are proportional in terms of their associated sets \(Z(A)\).

\begin{lemma}\label{lem:1}
Let $A\ne0$ and $C\in\B$. Then the following statements are equivalent: 
\begin{enumerate}
    \item The operators $A$ and $C$ are linearly independent, write $C=\lambda A$ for some $\lambda \in \C$.
    \item $Z(A)\subseteq Z(C)$.
    \item For every $R\in\Fone(X)$ 
\[
 (AR)^2=0\quad\Longrightarrow\quad(CR)^2=0.
\]
\end{enumerate}
Consequently, if $A,C\ne0$, then $Z(A)=Z(C)$ if and only if $A$ and $C$ are proportional.
\end{lemma}
\begin{proof}
The implications $(1)\Rightarrow(2)$ and $(2)\Rightarrow(3)$ are
immediate. It remains to prove that $(3)\Rightarrow(1)$.

Assume that, for every rank-one operator $T=x\otimes f$, where
$x\in X$ and $f\in X^*$, the condition $(AT)^2=0$ implies
$(CT)^2=0$. By the rank-one operator identity established above,
we obtain
\[
f(Ax)=0\quad\Longrightarrow\quad f(Cx)=0,
\qquad x\in X,\ f\in X^*.
\]
By the Hahn--Banach theorem, this implication yields
\(
\ker A\subseteq\ker C
\)
and, for every $x\in X$, the vectors $Ax$ and $Cx$ are
linearly dependent. A standard linear dependence argument shows that $A$ and $C$ are linearly dependent. Consequently,
there exists $\lambda\in\mathbb C$ such that $C=\lambda A$,
which completes the proof.
\end{proof}

\begin{lemma}\label{lem:rank-one-proportionnelle}
Let $M, N \in \mathcal{F}_1(X)$ be two rank-one operators. Then $M$ and $N$ are linearly dependent if and only if, for every rank-one idempotent $P \in \mathcal{I}_1(X)$, the following equivalence holds:
\[
(MP)^2 = 0 \quad \iff \quad (NP)^2 = 0.
\]
\end{lemma}
\begin{proof}
The direct implication is immediate. Conversely,
let $M = x \otimes f$ and $N = y \otimes g,$
for some $x, y \in X \setminus \{0\}$ and $f, g \in X^* \setminus \{0\}$. 

For every $z \in X$ and $h \in X^*$ satisfying $h(z) = 1$, the operator $R = z \otimes h$ is a rank-one idempotent. Since
\[
(MR)^2 = f(z)h(x)\,MR \quad \text{and} \quad (NR)^2 = g(z)h(y)\,NR.
\]
The hypothesis implies that
\begin{equation}\label{Eq:Rone}
    f(z)h(x) = 0 \quad \iff \quad g(z)h(y) = 0.
\end{equation}
First, we show that $f$ and $g$ are linearly dependent. Let $z \in \ker f \setminus \{0\}$. By the Hahn-Banach Theorem (by extending a linear functional defined on $\operatorname{span}\{x, y, z\}$), we can find a functional $h \in X^*$ such that $h(z) = 1$, $h(x) \neq 0$, and $h(y) \neq 0$. By equation \eqref{Eq:Rone}, we obtain $g(z) = 0$. Hence, $\ker f \subseteq \ker g$, which implies that $g = \lambda f$ for some scalar $\lambda \neq 0$.

Assume, by way of contradiction, that $x$ and $y$ are linearly independent. Choose $h \in X^*$ such that $h(x) = 0$ and $h(y) = 1$. We claim that there exists $z\in X$ satisfying $h(z) = 1$ and $g(z) \neq 0$. Indeed, we distinguish two  cases.  If $g$ and $h$ are linearly dependent, we may simply take $z = y$. Otherwise, if $g$ and $h$ are linearly independent, there exists 
$$z_1 \in \ker g \setminus \ker h\quad \text{ and } \quad z_2 \in \ker h \setminus \ker g,$$ such that $h(z_1) = 1$ and $g(z_2) = 1$; then the vector $z = z_1 + z_2$ satisfies $h(z) = g(z) = 1$. 

Finally, let consider the rank-one idempotent $P = z \otimes h$. A direct computation gives $(MP)^2 = f(z)h(x)MP = 0,$
whereas
$
(NP)^2 = g(z)h(y)NP = g(z)^2 y \otimes h \neq 0.
$
This contradicts the hypothesis that $(MP)^2 = 0 \iff (NP)^2 = 0$. Therefore, $x$ and $y$ must be linearly dependent. Consequently, the operators $M$ and $N$ are linearly dependent.
\end{proof}

For an operator $C\in \B$, let denote, 
\[
\Gamma(C):=\{R\in\mathcal{N}_1(X):(CR)^2=0\}.
\]
This set will play a key role in characterizing operators of the form
$aI+bN$, where $a,b\in\mathbb{C}$ and $N$ is a rank-one nilpotent operator.

\begin{lemma}\label{lem:gamma}
Assume that $\dim X \ge 2$, and let $N = u \otimes f \in \mathcal{N}_1(X)$, in particular $f(u) = 0$. For any $C \in \B$, the following equivalence holds:
\begin{equation}\label{eq:gamma}
\Gamma(C) = \Gamma(N) \quad \Longleftrightarrow \quad C = \alpha I + \beta N \quad \text{for some } \alpha \in \mathbb{C},\ \beta \in \mathbb{C}^*.
\end{equation}
\end{lemma}

\begin{proof}
First, suppose that $C=\alpha I+\beta N$, where $\alpha,\beta\in\mathbb{C}$ and
$\beta\ne0$. Let $M=x\otimes g\in\mathcal{N}_1(X)$, with $x\in X$ and
$g\in X^*$ satisfying $g(x)=0$. Then, for every operator $A\in \B$
\[
M\in\Gamma(A)
\quad\Longleftrightarrow\quad
(AM)^2=0
\quad\Longleftrightarrow\quad
g(Ax)=0.
\]
Since $g(Cx)=\beta g(Nx)$ and  $\beta\ne0$, it follows that $g(Cx)=0$ if and only if $g(Nx)=0$.
Consequently,
\[
M\in\Gamma(C)\quad\Longleftrightarrow\quad M\in\Gamma(N),
\]
and hence
\(
{\Gamma(C)=\Gamma(N).}
\)

For the converse, we first assume that $\dim X \ge 3$. Since $\ker f=\ker N$, we have $Nx=0$ for every $x\in\ker f$.
Let $x\in\ker f$ be a nonzero vector and let $g\in X^*$ be any nonzero linear
functional satisfying $g(x)=0$. Then the rank-one and  nilpotent operator
$R=x\otimes g$ satisfies $NR=0$, and hence $R\in\Gamma(N)$.
By our assumption, it follows that $R\in\Gamma(C)$, which implies $g(Cx)=0$ for every $g\in X^*$ such that $g(x)=0$. By the Hahn-Banach theorem, $Cx$ and $x$ are linearly dependent.
Since this holds for every $x \in \ker f$, the restriction of $C$ to $\ker f$ is a scalar multiple of the identity. Put $C\vert{}_{\ker f} = \alpha I$ for some scalar $\alpha$. In particular, the operator $T = C - \alpha I$ on $X$ is of rank at most one, and $\ker f \subset \ker T$. Consequently, we obtain $T=v \otimes f$ for some vector $v \in X$, and so $C = \alpha I + v \otimes f$.

For $x\notin\Ker f$, then for every $g\in X^*$ with $g(x)=0$, the operator $M=x\otimes g\in \None(X)$ and moreover, 
$$M\in \Gamma(C)\quad\Longleftrightarrow \quad M\in \Gamma(N).$$
On the other hand,  equation (\ref{eq:trace}) implies 
$$M\in \Gamma(C)\quad\Longleftrightarrow \quad \tr(MC)=0 \quad\Longleftrightarrow \quad g(v)f(x)=0,$$
Similarly,  $ M\in \Gamma(N)$ if and only if $g(u)f(x)=0.$ Since $x\notin \ker f$ and  $\Gamma(C)=\Gamma(N)$, it follows that, $g(v)=0$ if and only if $g(u)=0$. The Hahn-Banach theorem implies that $v \in \operatorname{span}\{x, u\}$ for every $x \notin \ker f$. Taking the intersection over these subspaces shows that $u$ and $v$ are linearly dependent, and thus  $v = \beta u$ for some scalar $\beta$.

To complete the proof of this part, we need to show that $\beta\ne 0$. Suppose  by contradiction that  $\beta = 0$. Then in particular $C=\alpha I$,  and the assumption $\Gamma(C) = \Gamma(N)$ implies that $\Gamma(\alpha I) = \mathcal{N}_1(X) = \Gamma(N)$. Consequently, for every $x\in X$ and $g\in X^*$,  the following implication holds 
$$g(x)=0\quad\Longrightarrow \quad g(u)f(x)=0.$$
Since $N = u \otimes f\in \None(X)$ is non zero,  there exists an element $x_0 \in X$ such that $f(x_0) \neq 0$. Hence, every continuous linear functional $g \in X^*$ that annihilates $x_0$ also annihilates $u$, again the Hahn-Banach Theorem, implies that $u=\lambda x_0$ for some scalar $\lambda\in \C$. Since $f(u)=0$ and $f(x_0)\ne 0$, then $\lambda=0$, whence $u=0$. This yields  $N=0$, which contradict the fact that $N\in \None(X)$. Therefore $\beta \ne 0$.

If $\dim X = 2$, identifying $\B$ with $\mathcal{M}_2(\mathbb{C})$. It is clear that if $X\in \mathcal{M}_2(\mathbb{C})$, then 
$$X^2=0\quad \iff\quad \tr(X)=0 ~\text{ and }~ det(X)=0 \quad \iff \quad X=\begin{pmatrix}x&y\\z&-x\end{pmatrix}\in M_2(\C),\quad x^2+yz=0.$$
Now, without loss of generality, we may assume that $N$ is of the form  
$$N=\begin{pmatrix} 0 & 1 \\ 0 & 0 \end{pmatrix},\quad\text{ and thus } \quad \Gamma(N) = \left\{ \begin{pmatrix} 0 & y \\ 0 & 0 \end{pmatrix} : y \in \mathbb{C}^* \right\} .$$
Let \(C=\begin{pmatrix}a&b\\c&d\end{pmatrix}\in M_2(\mathbb C)\).
A straightforward computation yields
 \[
\Gamma(C)=
\left\{
\begin{pmatrix}x&y\\z&-x\end{pmatrix}
\ne 0:
\begin{cases}
x^2+yz=0,\\
(a-d)x+bz+cy=0.
\end{cases}
\right\}
\]
Hence $\Gamma(C)=\Gamma(N)$ if and only if $a=d$, $c=0$ and $b\ne 0$. Equivalently, $C=aI+bN$, and this completes the proof of the result in the two-dimensional case.
\end{proof}

We denote by $\mathcal{U}_1$ the set of all operators in $\mathcal{B}(H)$ that can be written as a scalar multiple of the identity plus a rank-one nilpotent operator. We have 
\[
 \mathcal U_1=\{\alpha I+N:\alpha\in\C^*,\ N\in\None(X)\}.
\]
Every element of $\mathcal U_1$ is invertible, since
\begin{equation}\label{eq:inverseU}
 (\alpha I+N)^{-1}=\alpha^{-1}I-\alpha^{-2}N.
\end{equation}
The following lemma provides a characterization of rank-one operators in terms of the nilpotence of operator products and the set $\mathcal{U}_1$

\begin{lemma}\label{prop:rankonecriterion}
    Let $A\in\B$  be a nonzero operator on a Banach space $X$ with  $\dim X\ge2$. The following statements are equivalent:
\begin{enumerate}[label=(\roman*)]
 \item $\rank A=1$;
 \item $(AB)^2=0$ if and only if $(BA)^2=0$ for every $B\in\B$, and there exists $U\in\mathcal U_1$ such that $(AU)^2=0$.
\end{enumerate}
\end{lemma}
\begin{proof}
Suppose $A=x\otimes f\ne0$ is a rank-one operator.  For every $B\in \B$, we have 
$$ (BA)^2=0\iff f(Bx)=0 \iff (AB)^2=0,$$
and consequently, the first condition in (ii) holds. 

If $f(x)=0$, take $U=I+A$; then $A^2=0$, $U\in\mathcal U_1$, and $(AU)^2=0$. If $f(x)\ne0$, choose $y$ linearly independent of $x$ with $f(y)\ne0$, and choose $g\in X^*$ such that $g(y)=0$ and $g(x)=-f(x)/f(y)$. Then $N=y\otimes g$ is a rank-one nilpotent and $U=I+N\in\mathcal U_1$. Since $f(Ux)=f(x)+g(x)f(y)=0$, we have $(AU)^2=0$.

Conversely, suppose that (ii) holds. Since $U$ is invertible and $(AU)^2=0$, then $\Ker A\ne\{0\}$.  If $\rank A\ge2$, Choose $0\ne x_0\in\Ker A$, and $x_1,x_2\in X$ such that $y_1=Ax_1$ and $y_2=Ax_2$ are linearly independent. By Hahn-Banach Theorem, there exist  $f_1,f_2\in X^*$ with $f_i(y_j)=\delta_{ij}$. Set 
\[
 B=x_0\otimes f_1+x_1\otimes f_2.
\]
Then, $AB=y_1\otimes f_2$, so $(AB)^2=0$.  On the other hand  $BAx_2=x_1$ and $BAx_1=x_0$, and hence $(BA)^2x_2=x_0\ne0$. This contradicting our assumption in (ii). Thus $\rank A=1$.
\end{proof}
\begin{remark}
    The result of Lemma~\ref{prop:rankonecriterion} remains valid  when $\dim(X) = 2$. Indeed, since all elements of $\mathcal{U}_1$ are invertible, the condition $(AU)^2 = 0$ implies that $\operatorname{rank}(A) \le 1$ for every $U\in \mathcal U_1$. Since $A \neq 0$, it necessarily follows that $\operatorname{rank}(A) = 1$. Conversely, if $\operatorname{rank}(A) = 1$, the construction of $U$ in the proof of the preceding lemma does not rely on the assumption that $\dim(X) \ge 3$; hence, by the same arguments, we can find $U \in \mathcal{U}_1$ satisfying condition~(ii).
\end{remark}

\section{Auxiliary Results and Intermediate Lemmas}\label{section:lemmas}
In this section, we collect several auxiliary results concerning surjective maps $\Phi \colon \B \to \B$ that preserve the square-zero property of operator products in both directions. These results play a key role in the proof of our main theorems.

 \begin{lemma}\label{prop:projective}
If a surjective map $\Phi:\B\to\B$ satisfies \eqref{eq:P}, then:
\begin{enumerate}[label=(\roman*)]
 \item $\Phi(0)=0$, and $\Phi(A)\ne0$ whenever $A\ne0$;
 \item $A,C\ne0$ are proportional if and only if $\Phi(A),\Phi(C)$ are proportional;
 \item $A^4=0$ if and only if $\Phi(A)^4=0$.
\end{enumerate}
\end{lemma}

\begin{proof}
Condition \eqref{eq:P} and the surjectivity of $\Phi$ yield
\begin{equation}\label{Eq:Z1}
     Z(\Phi(A))=\Phi(Z(A))\qquad(A\in\B).
\end{equation}
Moreover, 
\begin{equation}\label{Eq:Z2}
    Z(A)=Z(C), \quad \text{if and only if} \quad Z(\Phi(A))=Z(\Phi(C)).
\end{equation}
Now $Z(0)=\B$, whereas $Z(A)\ne\B$ for $A\ne0$. It follows that $\Phi(0)=0$ and $\Phi(A)\ne 0$ for every nonzero operator  $A\in \B$.

Equation~\eqref{Eq:Z2} and Lemma~\ref{lem:1} show that $A$ and $C$ are proportional if and only if $Z(A) = Z(C)$, which  is equivalent to $Z(\Phi(A)) = Z(\Phi(C))$, and thus $\Phi(A)$ and $\Phi(C)$ are proportional.

Finally, taking $B=A$ in \eqref{eq:P} gives $A^4=0$ if and only if $\Phi(A)^4=0$.
\end{proof}

\begin{lemma}\label{cor:nilpreserved}
If  a surjective map $\Phi:\B\to\B$ satisfies \eqref{eq:P}, then for every $N\in \B$
$$N\in\None(X)\quad \iff \Phi(N)\in \None(X),$$
in particular, 
$$\Phi(\None(X))=\None(X).$$
\end{lemma}

\begin{proof}
Let $N \in \mathcal{B}(X)$ be a nonzero operator. By Lemma~\ref{lem:nilpotents},  the operator $N$ belongs in $\None(X)$ if and only if $N^4 = 0$ and 
$$(NA)^2 =0 \quad\iff\quad  (AN)^2 = 0, \quad\text{ for every } \quad A \in \B.$$ By  Lemma~\ref{prop:projective} and condition \eqref{eq:P}, this is equivalent to saying that $\Phi(N)^4 = 0$ and 
$$(\Phi(N)\Phi(A))^2 =0\quad\iff  \quad(\Phi(A)\Phi(N))^2 = 0, \text{ for every } \quad A \in \mathcal{B}(X).$$
Consequently, the surjectivity of $\Phi$ and  Lemma~\ref{lem:nilpotents} implies that $\Phi(N) \in \None(X)$.
\end{proof}

\begin{lemma}\label{lem:Upreserved}
Every  surjective map $\Phi:\B\to\B$  satisfies \eqref{eq:P} preserves  the set $\mathcal{ U}_1$ in both direction: for every $A\in \B$
$$A\in \mathcal{U}_1\quad \iff \quad \Phi(A)\in \mathcal U_1.$$
In particular, 
$$\Phi(\mathcal U_1)=\mathcal U_1.$$ 
\end{lemma}

\begin{proof}
By Corollary~\ref{cor:nilpreserved} and the condition (\ref{eq:P}), we deduce that for every $A\in \B$,  
$$\Phi(\Gamma(A))=\Gamma(\Phi(A)).$$
Moreover, for every $A, B\in \B$
$$\Gamma(A)=\Gamma(B)\quad\iff\quad \Gamma(\Phi(A))=\Gamma(\Phi(B)).$$
Let  $C\in \mathcal U_1$, and write $C=\alpha I+N$ for some $N\in\None(X)$ and nonzero scalar $\alpha \in \C^*$. Hence $C^4=\alpha^4I+4\alpha^3 N\ne 0$ and by Lemma~\ref{lem:gamma}, we have $\Gamma(C)=\Gamma(N)$.  It follows that $\Gamma(\Phi(C))=\Gamma(\Phi(N))$. Combining Lemma~\ref{lem:gamma} with the fact that $\Phi(N) \in \None(X)$ and $\Phi(C)^4\ne0$, we deduce that 
$$\Phi(C)=aI+b\Phi(N), \quad \text{ for some } \quad a,b\in \C\setminus\{0\}.$$ 
To prove the reverse inclusion, it suffices to follow the same arguments as in the first part.
\end{proof}

\begin{lemma}\label{thm:rankone}
Every surjective map $\Phi$ satisfying \eqref{eq:P} preserves the set of rank-one operators $\Fone(X)$  as well as the set of  rank-one non-nilpotent operators $\Fone(X)\setminus\None(X)$ in both directions:
\[A\in \Fone(X) \quad\Longleftrightarrow\quad\Phi(A)\in \Fone(X),
\]
and for every $A\in \Fone(X)$
$$A^2\ne 0\quad \iff \quad \Phi(A)^2\ne 0.$$
In particular, $$\Phi(\Fone(X))=\Fone(X)\quad \text{ and } \quad \Phi(\Fone(X)\setminus \None(X))=\Fone(X)\setminus \None(X).$$
\end{lemma}

\begin{proof}
By Lemma~\ref{prop:rankonecriterion}, an operator $A \in \B$ is of rank one if and only if, for every $B \in \B$, $(AB)^2 = 0 \iff (BA)^2 = 0$, and there exists $U \in \mathcal{U}_1$ such that $(AU)^2 = 0$. By Condition~\eqref{eq:P} and Lemma~\ref{lem:Upreserved}, this property is preserved by $\Phi$ in both directions. Therefore, $A$ is of rank one if and only if $\Phi(A)$ is as well. The last preserving property follows from the previous result and the fact that, for every $A \in \mathcal{F}_1(X)$, $$A^2 = 0 \iff A^4 = 0 \iff \Phi(A)^4 = 0 \iff \Phi(A)^2 = 0.$$
\end{proof}

\begin{lemma}\label{thm:rankone}
If a surjective map $\Phi$ satisfying \eqref{eq:P}, then for every $P\in \Fone(X)\setminus \None(X)$, we have that $tr(\Phi(P))\ne 0.$  Furthermore,
 the induced map $\psi :\Ione(X)\mapsto \Ione(X)$ defined by $$\psi(P)=\frac{1}{\tr (\Phi(P))}\Phi(P),$$ is  bijective. Moreover, $\psi$ preserves  the orthogonality between the idempotent in both directions:
$$PQ=QP=0\quad\Longleftrightarrow\quad\psi(P)\psi(Q)=\psi(Q)\psi(P)=0.$$
\end{lemma}
\begin{proof}
    For an arbitrary $P\in \Fone(X)$, we have $P\notin\None(X)$ if and only $P^2\ne 0$. By equation (\ref{eq:trace}), this is equivalent to $\tr(P)\ne 0$. The same equation~\eqref{eq:trace} implies that for every $P \in \mathcal{F}_1(X)$,   $P^2 = P$ if and only if $\operatorname{tr}(P) = 1$.
    
    Since $\Phi$ preserves the set of non-nilpotent rank-one operators in both directions, it follows that $\tr(\Phi(P))\ne0$, and then  $\tr(\psi(P))=1$ for every $P\in \Ione(X)$. Hence the maps $\psi$ sends every rank-one idempotent to a rank-one idempotent, and thus  $\psi$ is well defined. Moreover, the surjectivity of $\Phi$ implies that $\psi$ is also surjective. 
    
    To prove the injectivity of $\psi$, let $P, Q\in \Ione(X)$, and suppose that   $\psi(P)=\psi(Q)$. Hence, $Z(\Phi(P))=Z(\Phi(Q))$, and by 
    condition (\ref{eq:P}), this is equivalent to $Z(P)=Z(Q)$. By Lemma \ref{lem:1}, it follows that $P=tQ$ for some non zero scalar $t\in \C$, since $P,Q$ are both ideompotent, we must have $t=1$ and thus $P=Q$,  which establishes the injectivity, and hence the bijectivity of the maps $\psi$.  

    According to Lemmas~\ref{lem:adjacency} and~\ref{prop:projective} , together with condition~\eqref{eq:P}, the relation $\sim$ defined on the set of rank-one operators is preserved by $\Phi$ in both directions. Indeed, for two rank-one operators $R_1, R_2 \in \mathcal{F}_1(X)$, we have $R_1 \sim R_2$ if and only if there exists $R \in \mathcal{F}_1(X)$ such that $Z(R_1) \cap Z(R_2) \subset Z(R)$, which is equivalent to $Z(\Phi(R_1)) \cap Z(\Phi(R_2)) \subset Z(\Phi(R))$  and hence $\Phi(R_1) \sim \Phi(R_2)$. Therefore, the preservation of  orthogonality  follows directly from Lemma~\ref{lem:orthogonality}.
    \end{proof}
 Recall that an operator $P \in \mathcal{B}(X)$ is an idempotent if $P^2 = P$. Two idempotents $P, Q \in \mathcal{I}_1(X)$ are said to be orthogonal if $PQ = QP = 0$.
    In \cite{Semrl2005},  (see also \cite{Molnar2002b, Semrl2004}), the author provides a characterization of bijective maps on the set of rank-one idempotents that preserve orthogonality between idempotents in both directions. This result plays a crucial role in the proof of our main theorems. For completeness, we record this result in the following theorem:
\begin{theorem}[\cite{Semrl2005}]\label{thm:semrl}
Let $X$ be a Banach space with  $\dim X\ge3$. Let $\psi:\Ione(X)\longrightarrow\Ione(X)$ be a bijective maps preserves the orthogonality between  idempotents in both directions. Then $\Psi$ has one of the following forms.
\begin{enumerate}[label=(\roman*)]
 \item If $X$ is infinite-dimensional, then either
 \[
 \psi(P)=SPS^{-1},
 \]
 where $S:X\to X$ is bounded, invertible, and linear or conjugate-linear, or
 \[
 \psi(P)=SP^*S^{-1},
 \]
 where $P^*:X^*\to X^*$ is the Banach adjoint and $S:X^*\to X$ is a bounded invertible linear or conjugate-linear map. In this case $X$ must be reflexive.
 \item If $X$ is finite-dimensional, then
 \[
 \psi(P)=S P^\sigma S^{-1}\quad\text{or}\quad
 \psi(P)=S(P^\sigma)^tS^{-1},
 \]
 where $S$ is an invertible matrix,  $\sigma : \C\mapsto \C$ is an automorphism, $P^t$ denotes the transpose  of $P$ and $P^\sigma$ is the matrix whose entries are obtained by applying $\sigma$ to the entries of $P$.
\end{enumerate}
\end{theorem}

\section{The Proof of Main Theorems}\label{section:proof}

In this section, we complete the proofs of our main results, starting with the infinite-dimensional case. First, note that the  similarity transformations of the form $A \longmapsto S A S^{-1}$ preserve the square-zero property of operator products in both directions. Moreover, multiplying such a transformation by a nonzero scalar function $c \colon \mathcal{B}(X) \to \mathbb{C}^*$ yields a resulting transformation that still preserves this property. Therefore, we need only to prove the direct direction of Theorem~\ref{thm:infinite}

\begin{proof}[Proof of Theorem~\ref{thm:infinite}]
Assume first that $\Phi$ is a surjective map satisfying \eqref{eq:P}. By Lemma~\ref{thm:rankone} and Theorem~\ref{thm:semrl}, the restriction of $\Phi$ to the set of rank-one idempotents takes one of the following forms: there exists a function $\gamma : \mathcal{I}_1(X) \to \mathbb{C}^*$ such that either
\[
\Phi(P) = \gamma(P)SPS^{-1} \quad \text{for all } P \in \mathcal{I}_1(X),
\]
where $S : X \to X$ is a bounded, invertible, linear or conjugate-linear operator,  or either   $X$ is reflexive and 
\[
\Phi(P) = \gamma(P)SP^*S^{-1} \quad \text{for all } P \in \mathcal{I}_1(X),
\]
where $S : X^* \to X$ is a bounded, invertible, linear or conjugate-linear operator.

First, assume that $\Phi$ falls into the first case. Fix a rank-one operator $A \in \Fone(X)$, in particular $\Phi(A)$ is also of rank-one. By condition \eqref{eq:P}, for every rank-one idempotent $P \in \mathcal{I}_1(X)$, we have
\begin{align*}
(AP)^2 = 0 \quad &\Longleftrightarrow \quad (\Phi(A)\Phi(P))^2 = 0 \\
&\Longleftrightarrow \quad (\gamma(P)\Phi(A)SPS^{-1})^2 = 0 \\
&\Longleftrightarrow \quad (\Phi(A)SPS^{-1})^2 = 0 \\
&\Longleftrightarrow \quad ((S^{-1}\Phi(A)S)P)^2 = 0.
\end{align*}
By  Lemma~\ref{lem:rank-one-proportionnelle}, this is equivalent to saying that $A$ and $S^{-1}\Phi(A)S$ are linearly dependent. Since these operators are both of rank-one, there exists a nonzero scalar $c(A) \in \mathbb{C}^*$ such that $S^{-1}\Phi(A)S = c(A)A$, and hence
\[
\Phi(A) = c(A)SAS^{-1}, \qquad \text{for every $A \in \Fone(X)$.}
\]
To complete the proof of this part, it remains to show that the above form holds on the entire algebra $\mathcal{B}(X)$. To this end, combining Lemma~\ref{lem:1} with condition~\eqref{eq:P} and following the same arguments as before, we obtain the desired form on $\B$. 

We now show that $\Phi$ cannot take the second form. Suppose, to the contrary, that $\Phi$ is of the second form. By the same reasoning as before,  using the properties of the adjoint operation on Banach spaces, combined with Lemma~\ref{lem:rank-one-proportionnelle}, Lemma~\ref{lem:1}, and the reflexivity of $X$, we can extend this representation to the entire algebra $\B$.  Consequently, there exists a nonzero scalar-valued function $c:\B\to \C^*$ such that
\[
\Phi(A) = c(A)SA^*S^{-1}, \qquad \text{for every } A \in \mathcal{B}(X).
\]
Since \(\dim X\geq 3\), choose linearly independent vectors \(x_1,x_2,x_3\in X\) and nonzero functionals \(f_1,f_2,f_3\in X^*\) such that $f_i(x_j) = \delta_{i,j}$. Setting 
$$A = x_1 \otimes f_2 + x_2 \otimes f_3\quad \text{ and } \quad B = x_1 \otimes f_1 + x_2 \otimes f_2.$$ 
A straightforward computation yields
\[
AB = x_1 \otimes f_2, \quad \text{and thus} \quad (AB)^2 = (x_1 \otimes f_2)^2 = 0.
\]
On the other hand, we have $A^*B^* = (BA)^* = A^*$ and  $(A^*B^*)^2 = (A^*)^2=f_3\otimes x_1 \neq 0$. Consequently, $(\Phi(A)\Phi(B))^2 \neq 0$, whereas $(AB)^2 = 0$. This contradicts the condition \eqref{eq:P}, and hence the maps $\Phi$ cannot take this form. This completes the proof of 
Theorem~\ref{thm:infinite}.
\end{proof}

In the finite-dimensional setting, it is readily verified that every
transformation of the form
$A \mapsto c(A)S A^\sigma S^{-1}$, where $S$ is an invertible matrix
and $\sigma$ is an automorphism of $\mathbb{C}$, is well-defined and
satisfies our condition. Thus, it remains only to establish the direct
implication of Theorem~\ref{thm:finite}.
\begin{proof}[Proof of Theorem~\ref{thm:finite}]
First, suppose  that  $dim(X)=d\ge 3$. 
By Theorem~\ref{thm:semrl} and Lemma~\ref{lem:extendrankone}, there exist there exists a function $\gamma : \mathcal{I}_1(X) \to \mathbb{C}^*$ such that the restriction of $\Phi$ to the set of  rank-one idempotents takes one of the following two forms:
\[
P \longmapsto \gamma(P)SP^\sigma S^{-1} 
\quad \text{or} \quad 
P \longmapsto \gamma(P) S(P^\sigma)^t S^{-1}, 
\quad P\in \Ione(X).
\]
Using the same argument as in the infinite-dimensional case, the preceding forms can be extended to the set of all rank-one matrices, and $\Phi$ retains the same form on this set. To complete the proof of this part, we distinguish the following two cases:

\textbf{Case 1.} Suppose that $\Phi$ takes the first form on the set of rank-one matrices. For an arbitrary $A \in M_d(\mathbb{C})$, let 
\[
C = S^{-1}\Phi(A)S.
\]
Since $\Phi$ satisfies \eqref{eq:P}, applying this condition to $A$ and an arbitrary rank-one matrix $R$, yields
\[
(AR)^2 = 0 \quad \Longleftrightarrow \quad (\Phi(A)\Phi(R))^2=0 \quad \Longleftrightarrow \quad (CR^\sigma)^2 = 0.
\]
Applying the field automorphism $\sigma$ to the first equation, we obtain
\[
(A^\sigma R^\sigma)^2= 0 \quad \Longleftrightarrow \quad (CR^\sigma)^2 = 0.
\]
Since $R^\sigma$ ranges over the entire set of rank-one matrices as $R$ does, it follows from the Lemma~\ref{lem:1}  that  $C$ and $A^\sigma$ are proportional. Thus, there exists a nonzero scalar $c(A)$ such that $C = c(A)A^\sigma.$  Consequently,
\[
\Phi(A) = c(A)SA^\sigma S^{-1}.
\]

\textbf{Case 2.} Suppose now that $\Phi$ takes the second form on the set of rank-one matrices. The same argument, gives
\[
(AR)^2 = 0 \quad \Longleftrightarrow \quad \left(C(R^\sigma)^t\right)^2 = 0.
\]
Since for every rank-one matrix $R$ and every $C\in M_n(\C)$
\[
\left(C(R^\sigma)^t\right)^2=0 \quad \Longleftrightarrow \quad \left(R^\sigma C^t\right)^2=0\quad \Longleftrightarrow \quad \left(C^t R^\sigma\right)^2=0.
\]
and applying $\sigma$ to the first trace expression, Lemma~\ref{lem:1} yields
\[
C^t = c(A)A^\sigma,
\]
for some nonzero scalar $c(A)$. Equivalently,
\[
\Phi(A) = c(A)S(A^\sigma)^t S^{-1}.
\]
We show that this second possibility is impossible when $d \ge 3$. Indeed, consider the matrices $A$ and $B$ constructed in the proof of Theorem~\ref{thm:infinite}, which satisfy
\[
(AB)^2 = 0 \quad \text{and} \quad (BA)^2 \neq 0.
\]
Since $\sigma$ preserves multiplication and the transpose reverses the order of multiplication, we have
\[
\big((A^\sigma)^t(B^\sigma)^t\big)^2 = \big((B^\sigma A^\sigma)^2\big)^t \neq 0.
\]
This contradicts \eqref{eq:P}. Therefore, $\Phi$ cannot take  the second  expression.
\end{proof}

\end{document}